\documentclass{amsart}

\title[Multiplicities on Schubert Varieties]{Multiplicities of the Most Singular Point on Schubert Varieties on Gl(n)/B for $n=5,6$}
\author{Alexander Woo}
\begin{document}

\begin{abstract}
We calculate using Macaulay 2 the multiplicities of the most singular
point on Schubert varieties on Gl(n)/B for $n=5,6$.  The method of
computation is described and tables of the results are included.
\end{abstract}

\maketitle

\section{Introduction}

In this paper, we compute, using Macaulay 2, the multiplicity of the
most singular point on Schubert varieties of the flag manifold
$Gl(n)/B$ for $n=5$ and $n=6$, using a description of the preimage of
Schubert varieties in $GL(n)$ first given by Fulton \cite{Ful2} and
more recently developed by Knutson and Miller \cite{K-M}.  Results of
Krattenthaler \cite{Kra}, Rosenthal and Zelevinsky \cite{R-Z}, and
Lakshmibai and Weyman \cite{L-W} give combinatorial and determinantal
formulas for multiplicities (at all points) of Schubert varieties on
the Grassmanian.  Furthermore, the singular loci of Schubert varieties
on flag manifolds have been much studied, with known results collected
in \cite{B-L}, starting with the fundamental result of Lakshmibai and
Sandhya \cite{L-S} that a Schubert variety indexed by the permutation
$w$ is singular iff $w$ contains either the pattern $1324$ or the
pattern $2143$.  However, no results for the multiplicities of the
Schubert varieties appear to be known on the full flag variety.

In the next section, we briefly define the objects under study and
outline some basic results about them.  A detailed introduction to
Schubert varieties can be found, for example, in \cite[part III]{Ful}.
This section also serves to fix our conventions; the differing choices
made by different authors in the subject can cause significant
confusion.  In particular, our convention for indexing Schubert
varieties is opposite to the convention used in \cite{B-L}, which now
seems to be fairly standard.  The third section describes our
algorithm, and the fourth section demonstrates this algorithm for
$w=1324$.  Two appendices give our code and the results of our
computations.

\section{Definitions and Conventions}

A {\it (complete) flag} $\mathcal{F}$ in $\mathbb{C}^n$ is a sequence
of subspaces $\{0\}=F_0\subset F_1\subset F_2\subset\cdots\subset
F_n=\mathbb{C}^n$ such that the subspace $F_i$ has dimension $i$.
Fixing a basis for $\mathbb{C}^n$, we can represent $\mathcal{F}$ by
a matrix $M$ as follows.  For each component $F_i$ of $\mathcal{F}$,
pick a vector $m_i\in F_i\setminus F_{i-1}$, and write $m_i$ as the
$i$-th row of $M$.  Note that $F_i$ will be the span of the first $i$
rows of $M$; in particular, $M$ is in invertible since
$F_n=\mathbb{C}^n$.

This representation of the flag is clearly not unique, as it involves
repeated choices of vectors.  To be precise, both multiplying any row
of $M$ by a constant and adding any row to a subsequent row of $M$
leaves the flag unchanged.  This is equivalent to multiplying $M$ on
the left by a lower triangular matrix.  We can now give the structure
of an algebraic variety to the set of flags; namely, it is the
quotient of $G=Gl(n)$ by the group of lower triangular
matrices $B^-$ acting on the left, and we denote this variety by
$B^-\backslash G$ and call it the {\it flag variety}.  We also have a natural
map

$$\pi: G\rightarrow B^-\backslash G$$

sending a matrix to the flag it represents.

We can, however, pick a standard representation $N(\mathcal{F})$ for
each flag $\mathcal{F}$ as follows.  Pick a matrix $M$ that represents
$\mathcal{F}$.  Now take the leftmost nonzero entry in the first row
of $M=\left[m_{ij}\right]$, which we call $m_{1l}$, and use elements of
$B^-$ to change $M$ so that $m_{kl}=0$ for $k>1$.  We then take the
leftmost nonzero entry of the second row, place $0$s in all entries
below it, and repeat for all rows in order to finally get
$N(\mathcal{F})$.  Unfortunately, we cannot identify $B^-\backslash G$ with
these representations of flags except as a set, since this process of
row reduction destroys the geometry of $B^-\backslash G$.

Let the group $B^+$ of upper triangular matrices act on $G$ by right
multiplication; by associativity of multiplication, this commutes with
the left action of $B^-$ and therefore gives an action on $B^-\backslash G$.
This action can be thought of as adding any column of a matrix to a
column to its right.  Under these two actions, any matrix can be sent
to a unique permutation matrix $W$, so we can index orbits of the
right action of $B^+$ on $B^-\backslash G$ by permutations.  In particular, a
flag $\mathcal{F}$ is in the orbit of $\pi(W)$, where $W$ is the
permutation matrix with $1$s in the leftmost nonzero entries of each
row in $N(\mathcal{F})$.  Let $X_w$ denote the orbit of the flag
$\pi(W^{-1})$, where $W$ is the permutation $w$ written as a matrix.
$X_w$ is known as a {\it Schubert cell}.  The process of choosing the
representative $N(\mathcal{F})$ for each flag does preserve geometry
locally on $X_w$, so $X_w$ is isomorphic to
$\mathbb{A}^{\binom{n}{2}-l(w)}$, where $l(w)$ is the length of a
shortest expression of $w$ as a product of adjacent transpositions,
or, equivalently, the number of inversions in $w$.

We denote by $Y_w$ the closure of $X_w$ in $B^-\backslash G$; it is known as a
{\it Schubert variety}.  For permutations $v,w\in S_n$, let $v\succ w$
if $l(v)>l(w)$ and $v=tw$ for some transposition $t$.  The transitive
closure of the relation $\succ$ is known as the {\em Bruhat order};
for the remainder of this paper, $v>w$ for $v,w\in S_n$ means that $v$
is greater than $w$ in this partial order.  It is a classical result
that

$$Y_w = \bigcup_{w^\prime>w} X_{w^\prime}.$$

Note that the unique 0-dimensional cell $X_{w_0}$, where
$w_0=n\cdots21$ is contained in every Schubert variety.

Given a variety $X$, the {\it multiplicity} at a point $p$ of $X$,
which we will denote $\operatorname*{mult}_p X$ is the degree of the
projective tangent cone
$\operatorname{Proj}(\operatorname*{gr}_{\mathfrak{m}_p}
\mathcal{O}_{X,p})$ as a subvariety of the projective tangent space
$\operatorname{Proj} (\operatorname*{Sym}^*
\mathfrak{m}_p/\mathfrak{m}^2_p)$, or, equivalently, if the
Hilbert--Samuels polynomial of $\mathcal{O}_{X,p}$ is written
$a_nx^n+a_{n-1}x^{n-1}+\cdots+a_0$, $\operatorname*{mult}_p X=n!a_n$.  For a
nonsingular point $p$, $\operatorname*{mult}_p X=1$; at a singular point $q$,
$\operatorname*{mult}_q X>1$ and measures roughly how singular $X$ is at $q$.
Slightly more precisely, the multiplicity counts how many times a
generic hyperplane cuts through $X$ in a neighborhood of $q$.

Since a Schubert variety is invariant under the right action of $B^+$,
its multiplicity must remain constant on $B^+$ orbits, or Schubert
cells.  Moreover, since for $v\succ w$, there exists a $\mathbb{P}^1$
with one point in $X_v$ and the remaining points in $X_w$, by
semicontinuity, multiplicities must be nondecreasing with respect to
Bruhat order.  $X_{w_0}$ must therefore be the most singular
point of $Y_w$ and the multiplicity of $Y_w$ there measures how
singular $Y_w$ gets.  In particular, a multiplicity of $1$ at
$X_{w_0}$ indicates that $Y_w$ is smooth.

\section{Explanation of the Algorithm}

Since multiplicity is a local property, we can calculate it after
restricting to an affine neighborhood of $X_{w_0}$ in $B^-\backslash G$.  A
natural candidate is $\Omega_{w_0}$, the orbit of $W_0$ (defined to be
$w_0$ written as a permutation matrix and considered as a flag) under
the right action of $B^-$.  (In general, $\Omega_w$ is defined as the
orbit of the flag $\pi(W^{-1})$ under the right action of $B^-$
(rather than $B^+$) and is known as a {\it dual Schubert cell}.)
Locally on $\Omega_{w_0}$, the map $\pi: G\rightarrow B^-\backslash G$ has a
section $\sigma$, namely the map that sends a flag $\mathcal{F}$ to
$N(\mathcal{F})$.  This identifies $\Omega_{w_0}$ with the matrices
with $1$s on the main antidiagonal and $0$s to the right and below;
$X_{w_0}$ is mapped to the permutation matrix $W_0$.  Since $\sigma$
is a local section, $Y_w$ is locally isomorphic in a neighborhood of
$X_{w_0}$ to $\pi^{-1}(Y_w)\cap\sigma(\Omega_{w_0})$.

Now we need to find equations defining $\pi^{-1}(Y_w)$.  Fix a
permutation $w\in S_n$.  Let $R(w)=\left[r_{ij}(w)\right]$ be the
integer matrix with $r_{ij}(w)=\#\{w^{-1}(k)\leq i,k\leq j\}$.  For any
matrix $M$, let $M_{ij}$ denote the submatrix consisting of the first
$i$ rows and first $j$ columns.  Then, for any invertible matrix $M$
with $\pi(M)\in X_w$, the rank of the submatrix $M_{ij}$ will be
$r_{ij}$.  The proof of this claim is as follows.  Note that the rank of
$M_{ij}$ for any $i$ and $j$ does not change under multiplication by
$B^-$ on the left, since the effect of multiplication by $b\in B^-$ on
the first $i$ rows is the same as that of multiplying by an element of
$Gl(i)$, namely the submatrix of $b$ consisting of the first $i$ rows
and columns.  Therefore, the claim can be verified on matrices of the
form $N(\mathcal{F})$ for $\mathcal{F}\in X_w$, where it is trivial.

It is a nontrivial combinatorial fact that Bruhat order can be
equivalently defined by $v>w$ if $r_{ij}(v)\leq r_{ij}(w)$ for all $i,j$.
Therefore $\pi^{-1}(Y_w)$ consists of all invertible matrices $M$
satisfying the rank conditions $\operatorname{rk}(M_{ij})\leq r_{ij}(w)$.
Now let $Z=\left[z_{ij}\right]$ be a matrix of indeterminates, and
$I_w$ be the ideal generated by all size $1+r_{ij}(w)$ minors of
$Z_{ij}$, for all $i$ and $j$.  (As shown in \cite{Ful2}, where $I_w$
was originally defined, a small subset of the minors suffices to
generate $I_w$, but we will not use this fact.)  By the above
statement, it is clear that $I_w$ vanishes precisely on the points of
$\pi^{-1}(Y_w)$; Knutson and Miller \cite{K-M} prove that $I_w$ is in
fact radical, so $I_w=I(\pi^{-1}(Y_w))$.  Now, by sending $z_{ij}$ to
$0$ for entries below the main antidiagonal and $1$ for entries on the
main antidiagonal, we have the ideal $J_w$ for
$\pi^{-1}(Y_w)\cap\sigma(\Omega_{w_0})$ as a subvariety of
$\sigma(\Omega_{w_0})\cong \mathbb{A}^{\binom{n}{2}}$.

In our coordinates for $\sigma(\Omega_{w_0})$, $X_{w_0}$ corresponds
to the point where $z_{ij}=0$ for all $i$ and $j$, or equivalently the
maximal graded ideal $\mathfrak{m}=\langle z_{ij}\rangle$.  The
projective tangent cone of $Y_w$ is therefore $\operatorname{Proj}
(\operatorname*{gr}_\mathfrak{m} S/J_w)$, where $S=k[z_{ij}]$.  For a
polynomial $f\in S$, let $f_d$ be the homogeneous part of degree $d$,
and let $s(f)$ be the smallest number such that $f_{s(f)}\neq 0$.  Let
$J^\prime_w=\langle f_{s(f)} | f\in J_w\rangle$; then
$\operatorname*{gr}_m S/J_w\cong S/J^\prime_w$, since any polynomial
$f$ is sent to $f \pmod{\mathfrak{m}^{s+1}}$, of which $f_s$ is a
representative.  Degree is invariant under Gr\"obner deformation and
easiest to calculate on monomial ideals, so, to calculate the
multiplicty of $Y_w$, or, equivalently, the degree of $S/J^\prime_w$,
it suffices to calculate the degree of
$S/\operatorname{in}(J^\prime_w)$ under any (graded) term order.  This
is the same as taking the initial ideal of $J_w$, taking care to use a
term order which places the {\bf lowest} degree term first.

In practice, on Macaulay 2, one does this by homogenizing the
generators of $J_w$ using a new variable $t$ (or by sending $z_{ij}$
to $t$ rather than $1$ for entries on the main diagonal in passing
from $I_w$ to $J_w$), and using a term order that refines the partial
order by degree in $t$.  (See \cite[Prop 15.28]{Eis} for a proof that
this is equivalent.)  We can then compute the initial ideal, send $t$
to $1$, and calculate the degree.

\section{Example}

As an example, we calculate the multiplicity of $Y_w$ at $X_{w_0}$ for
$w=2143$, the smallest nontrivial example.  Here, we have

$$W^{-1}=\left[
\begin{matrix}
0&1&0&0\\
1&0&0&0\\
0&0&0&1\\
0&0&1&0\end{matrix}\right].$$

Therefore, we have the rank matrix

$$R=\left[
\begin{matrix}
0&1&1&1\\
1&2&2&2\\
1&2&2&3\\
1&2&3&4\end{matrix}\right].$$

Exactly two of the rank entries give minors in $I_w$, namely $r_{11}=0$,
which gives $z_{11}\in I_w$, and $r_{33}=2$, which gives

$$\det\left[\begin{matrix}
z_{11}&z_{12}&z_{13}\\
z_{21}&z_{22}&z_{23}\\
z_{31}&z_{32}&z_{33}\end{matrix}\right]=0.$$

Setting $z_{33}=0$ and $z_{32}=z_{23}=1$, we have $J_w=\langle z_{11},
-z_{11}+z_{12}z_{31}+z_{21}z_{13}-z_{31}z_{22}z_{13}\rangle=\langle
z_{11}, z_{12}z_{31}+z_{21}z_{13}-z_{31}z_{22}z_{13}\rangle$.  Then
$J^\prime_w=\langle z_{11}, z_{12}z_{31}+z_{21}z_{13}\rangle$, and it
is clear that the multiplicity of $Y_w$ at $X_{w_0}$ is $2$.  For
purposes of illustration, we carry out the remainder of the algorithm.
Homogenizing the generators of $J_w$ gives us $\langle z_{11},
tz_{12}z_{31}+tz_{21}z_{13}-z_{31}z_{22}z_{13}\rangle$, and one
possible appropriate initial ideal is $\langle z_{11},
tz_{21}z_{13}\rangle$.  Sending $t$ to $1$, we get $\langle z_{11},
z_{21}z_{13}\rangle$, which has degree 2; therefore, we conclude that 
$\operatorname*{mult}_{X_{w_0}} Y_w = 2.$

\section{Acknowledgements}

This paper was originally a project for Bernd Sturmfel's Math 274
class in Fall 2002.  Thanks also go to Ezra Miller for suggesting this
project and helping me get started.

\appendix
\section{Maple and Macaulay 2 code}

The following is our Maple code for generating the rank matrices.

\begin{verbatim}
with(combinat);
with(linalg);
interface(prettyprint=false);
n:=5;

sum(matrix(n,n,
      (i,j)->`if`((j>=k) and (i>=op(k, pPerm)), 1, 0)), k=1..n);

\end{verbatim}

The Maple output was sent to a file, and all square brackets [] were
converted to curly brackets \{\} for Macaulay 2.  Macaulay 2 code is
as follows.  Note that the first five lines must be changed for each $n$.

\begin{verbatim}

R=QQ[t,x11,x12,x21,x13,x22,x31,x14,x23,x32,x41,
       MonomialOrder=>Eliminate 1];

G = matrix({{x11,x12,x13,x14,t},
            {x21,x22,x23,t,0},
            {x31,x32,t,0,0},
            {x41,t,0,0,0},
            {t,0,0,0,0}})

S=QQ[x11,x12,x21,x13,x22,x31,x14,x23,x32,x41];

f=map(S,R,{1,x11,x12,x21,x13,x22,x31,x14,x23,x32,x41});

n=5;

# Mlist = (paste from Maple output)

# compute J_w
Ilist = apply(Mlist,
          M->trim(sum(flatten(for i from 0 to n-1 list
                                for j from 0 to n-1 list
                                  minors(M_(i,j)+1, 
                                    submatrix(G, {0..i}, {0..j})))));

GBlist = apply(Ilist, gb);

# gives in(J_w)
LTlist = apply(GBlist, GB -> leadTerm(gens(GB)));

# gives in(J^\prime_w)
ELTlist = apply(LTlist, LT->f(LT));

# outputs degrees
Dlist = apply(ELTlist, LT -> degree(ideal(LT))

\end{verbatim}

The output of the last line gives the degrees.  Output from other
lines, such as the initial ideals, could also be of interest.  These
computations run quite quickly; in fact, generating the ideals $J_w$
was by far the slowest step, and the Gr\"obner basis computation took
only a few seconds for $n=6$.

\newpage

\section{Computational Results for $n=5$ and $n=6$}

We have listed the permutations by multiplicity of $Y_w$ at $X_{w_0}$.
For each multiplicity, permutations are listed in lexicographic order.

\vspace{0.0in}

\begin{tabular}{c | l}
Multiplicity & Permutations\\
\hline
5 & 14325 \\
3 & 13425, 14235, 24153, 31524 \\
2 & 12435, 13245, 13254, 13524, 14253, 14352, 15324, 21354, 21435, \\
  & 21453, 21534, 21543, 23154, 24135, 24315, 24351, 25143, 31254, \\
  & 31425, 31542, 32154, 32514, 32541, 41325, 42153, 51324, 52143 \\
1 & 12345, 12354, 12453, 12534, 12543, 13452, 13542, 14523, 14532, \\
  & 15234, 15243, 15342, 15423, 15432, 21345, 23145, 23415, 23451, \\
  & 23514, 23541, 24513, 24531, 25134, 25314, 25341, 25413, 25431, \\
  & 31245, 31452, 32145, 32415, 32451, 34125, 34152, 34215, 34251, \\
  & 34512, 34521, 35124, 35142, 35214, 35241, 35412, 35421, 41235, \\
  & 41253, 41352, 41523, 41532, 42135, 42315, 42351, 42513, 42531, \\
  & 43125, 43152, 43215, 43251, 43512, 43521, 45123, 45132, 45213, \\
  & 45231, 45312, 45321, 51234, 51243, 51342, 51423, 51432, 52134, \\
  & 52314, 52341, 52413, 52431, 53124, 53142, 53214, 53241, 53412, \\
  & 53421, 54123, 54132, 54213, 54231, 54312, 54321

\end{tabular}

\vspace{0.03in}

For $n=6$, we have omitted the permutations of multiplicities 1 and 2
for reasons of space.  (There are 366 nonsingular Schubert varieties
and 207 with multiplicity $2$ at $X_{w_0}$.)  These can be inferred
from the remainder with the help of the pattern avoidance criterion.

\vspace{0.03in}

\begin{tabular}{c | l}
Multiplicity & Permutations\\
\hline
14 & 154326 \\
10 & 153426 \\
9 & 145326, 154236 \\
8 & 321654 \\
7 & 135426, 143526, 152436, 153246, 254163, 416325 \\
6 & 145236, 132546, 214365 \\
5 & 125436, 135246, 142536, 143256, 143265, 143625, 146325, 153264, \\
  & 154263, 154362, 164325, 215436, 251364, 251436, 251463, 253164, \\
  & 254136, 254316, 254361, 314625, 315426, 316425, 413625, 415326, \\
  & 426153, 514326, 614325 \\
4 & 153624, 152346, 134526, 214635, 215364, 215463, 216435, 231564, \\
  & 231654, 241365, 243165, 245163, 312645, 312654, 314265, 321564, \\
  & 321645, 326154, 351426, 351624, 413265, 416235, 421653 \\
3 & 124356, 124365, 124536, 125346, 132564, 132645, 132654, 134256, \\
  & 134265, 134625, 135264, 136425, 142356, 142365, 142635, 145263, \\
  & 145362, 146235, 152364, 152463, 153462, 163425, 164235, 214356, \\
  & 214536, 215346, 216453, 216543, 231546, 235164, 241536, 241563, \\
  & 241635, 241653, 245136, 245316, 245361, 246153, 251346, 251634, \\
  & 251643, 253146, 253416, 253461, 264153, 312546, 314526, 315246, \\
  & 315264, 315624, 316245, 316254, 316524, 321546, 325164, 341625, \\
  & 351264, 351642, 352164, 352614, 352641, 361524, 425163, 412635, \\
  & 413526, 415236, 416253, 416352, 421635, 423165, 426135, 426315, \\
  & 426351, 431625, 432165, 513426, 514236, 524163, 531624, 613425, \\
  & 614235, 624153, 631524

\end{tabular}

\end{document}